\documentclass{elsart}

\usepackage{graphicx}
\usepackage{amsfonts}
\usepackage{hyperref}
\usepackage{textcomp}

\makeatletter
\def\url@leostyle{%
  \@ifundefined{selectfont}{\def\UrlFont{\sf}}{%
    \def\UrlFont{\footnotesize\ttfamily}}\Url@do
}
\makeatother

\newtheorem{lemma}{Lemma}
\newtheorem{theorem}[lemma]{Theorem}

\newtheorem{observation}[lemma]{Observation}

\theoremstyle{definition}
\newtheorem{definition}{Definition}
\newtheorem{conjecture}[definition]{Conjecture}

\theoremstyle{remark}
\newtheorem{question}{Question}

\newcommand{\Aut}{\mbox{\sc Aut}}

\newcommand{\slot}{\Diamond}
\newcommand{\sslot}{\, \slot \,}
\newcommand{\M}{{\mathbf m}}
\renewcommand{\L}{{\mathbf l}}
\newcommand{\R}{{\mathbf r}}
\newcommand{\F}{{\mathbf f}}

\newcommand{\ind}{{\mathsf{ind}}}
\newcommand{\ws}[1]{L(#1)}

\newcommand{\B}{{\mathcal B}}
\newcommand{\Lock}{{\mathcal L}}

\renewcommand{\S}{{\mathcal S}}

\newcommand{\avoid}{\mathrm{Av}}
\newcommand{\intersect}{\cap}

\date{}

\begin{document}
\begin{frontmatter}

\title{On the Wilf-Stanley limit of $4231$-avoiding permutations and a conjecture of Arratia.}

\author{M.~H.~Albert}
\address{Department of Computer Science \\%
University of Otago, Dunedin, New Zealand}
\ead{malbert@cs.otago.ac.nz}
\author{M.~Elder\thanksref{EPSRC}}
\address{School of Mathematics and Statistics\\
University of St. Andrews, St. Andrews, Scotland} 
\ead{murray@mcs.st-and.ac.uk}
\thanks[EPSRC]{Supported by EPSRC grant GR/S53503/01.}
\author{A.~Rechnitzer\thanksref{ARC}}
\address{Department of Mathematics and Statistics\\
University of Melbourne, Melbourne, Australia} 
\ead{andrewr@ms.unimelb.edu.au}
\thanks[ARC]{The financial support of the Australian Research Council is gratefully acknowledged.}
\author{P.~Westcott}
\address{Melbourne, Australia}
\ead{p.westcott@gmail.com}
\author{M.~Zabrocki}
\address{Department of Mathematics and Statistics\\
York University, Toronto, Canada} 
\ead{zabrocki@mathstat.yorku.ca}

\begin{abstract}
We construct a sequence of finite automata that accept subclasses of
the class of $4231$-avoiding permutations. We thereby show that the
Wilf-Stanley limit for the class of $4231$-avoiding permutations is bounded
below by $9.35$. This bound shows that this class has the
largest such limit among all classes of permutations avoiding a
single permutation of length $4$ and refutes the conjecture that the
Wilf-Stanley limit of a class of permutations avoiding a single
permutation of length $k$ cannot exceed $(k-1)^2$.
\end{abstract}

\end{frontmatter}

\maketitle

\section{Introduction}

Let $\sigma = \sigma_1 \sigma_2 \cdots \sigma_k$ and $\pi = \pi_1
\pi_2 \cdots \pi_n$ be permutations of $\{1,2,3,\ldots,k \}$ and
$\{1,2,3,\ldots, n\}$ respectively, written as their sequences of
values. Then {\em $\sigma$ occurs as a pattern in $\pi$\/} if for some
subsequence $\tau$ of $\pi$ of the same length as $\sigma$ all the values in $\tau$ occur in the
same relative order as the corresponding values in $\sigma$. If
$\sigma$ does not occur as a pattern in $\pi$ we say that $\pi$ {\em
avoids\/} $\sigma$. A {\em pattern class of permutations}, or simply
class, is any set of permutations of the form:
\[
\avoid (X) = \{ \pi \, : \, \forall \, \sigma \in X, \: \mbox{$\pi$
avoids $\sigma$} \}
\]
where $X$ is any set of permutations. We usually write
$\avoid(\sigma)$ rather than $\avoid(\{ \sigma \})$. Pattern classes are the lower ideals of the set of all finite permutations with respect to the partial order ``occurs as a
pattern in'' and so are closed under arbitrary intersections and
unions.

Much of the study of pattern classes has concentrated on enumerating
classes $\avoid(X)$ when $X$ is a relatively small set of relatively
short permutations. We write $s_n(X)$ for $\left| \avoid(X) \cap {\S}_n \right|$. 
Results in this area led to the proposal of the
Wilf-Stanley conjecture. A somewhat simplified version of this
conjecture is:

\begin{conjecture}[Wilf-Stanley] Let $X$ be any non-empty set of
permutations. Then there exists a real number $c_X$ such that $s_n(X) \leq c_X^n$.
\end{conjecture}

The resolution of the Wilf-Stanley conjecture by Marcus and Tardos
\cite{MT:wilfStanley}, together with a result of Arratia's on the
classes defined by avoiding a single permutation
\cite{Arratia:growthRates} implies that for each permutation $\pi$
there exists a positive real number $\ws{\pi}$ called the {\em Wilf-Stanley limit\/} of the class $\avoid(\pi)$ such that:
\[
\lim_{n \to \infty} {s_n(\pi)}^{1/n} = \ws{\pi}.
\]
The values of $\ws{\pi}$ are known exactly for all
permutations of length $3$, and for all permutations of length $4$
except $4231$ and $1324$ (which have the same Wilf-Stanley limit, by the
obvious isomorphism between the corresponding classes). Using a result of Regev \cite{Regev:AdvMath1981}, B{\'o}na
\cite{Bona:JCTA97,Bona:DM97,Bona:AAM2004} provided bounds:
\[
9 \leq \ws{4231} \leq 288.
\]
Further results of B{\'o}na \cite{Bona:JCTA2005} show that $\ws{\pi} \geq (k-1)^2$ for all  layered permutations $\pi$ of length $k$ (a permutation $\pi$ is layered if $\pi_{j+1} < \pi_j$ implies $\pi_{j+1} = \pi_{j}-1$).
Arratia \cite{Arratia:growthRates} conjectured that for
all permutations $\pi$ of length $k$, $\ws{\pi} \leq (k-1)^2$. Regev's result shows that the
value $(k-1)^2$ is attained for the permutation $\pi = 123 \cdots
k$. In this paper we refute this conjecture by proving that:
\[
\ws{4231} \geq 9.35.
\]

Our proof of this result makes use of the insertion encoding
for permutations. We establish that there is a class of
permutations, strictly contained in $\avoid(4231)$ whose elements
are in one to one correspondence with the words of a language
accepted by a certain finite automaton. Using the standard transfer
matrix approach we are able to determine the growth rate of this
language, which thus provides the lower bound cited above.

In order to make this paper self-contained we provide a brief
introduction to the insertion encoding in the next section. Then we
will describe the automaton (actually a sequence of automata)
referred to above, and prove the required correspondence. We include
a brief discussion of the computational methodology and then a
summary and conclusions.

\section{The Insertion Encoding}

The insertion encoding is a general method for describing
permutations. It shares some similarity with the generating
tree approach of West \cite{We01,We02}  and also the enumeration schemes of Zeilberger \cite{Zeilberger:AnnComb1998} two approaches which have been used in these papers and elsewhere
 \cite{CM:AAM2002,EP:JCTA2004,Krattenthaler:AAM2001,KS:DM2003,SW:JAlgCom2002,SW:DM2004} to
enumerate  or determine structural information about a number of permutation classes.

A permutation $\pi$ is viewed as ``evolving'' by the successive
insertion of new maximal elements. Thus, the stages in the
evolution of $264153$ are: $\epsilon$ (the empty word), $1$, $21$, $213$,
$2413$, $24153$ and $264153$. Each step of the evolution is
described by a code letter of the form $\F_i$, $\L_i$, $\R_i$ or
$\M_i$ where $i$ is a positive integer. The intent of the symbols
will become more clear if in the evolution of $\pi$ we also include
placeholders, called {\em slots\/} in positions where an element will
eventually be inserted. We denote a slot by the symbol $\slot$. Now
the evolution of $264153$ can be written as:
\[
\slot  \to   \slot 1 \slot  \to  2 \slot 1 \slot  \to  2 \slot 1
\slot 3  \to  24 \slot 1 \slot 3  \to  24 \slot 153  \to 246153.
\]

It can be seen that each event in the evolution is of one of four
types: filling a slot (the last two events), insertion on the left
hand end of a slot (the addition of $4$), on the right hand end of a
slot (the addition of $3$), or in the middle of a slot splitting it
in two (the addition of $1$). The code letters then describe the
type of insertion to carry out, and the subscript denotes the slot
in which to perform the insertion (counted from left to right). Thus
the insertion encoding of $246153$ is $\M_1 \L_1 \R_2 \L_1 \F_2
\F_1$.

In considering $\avoid(4231)$ it turns out that a small modification
of this encoding provides a more natural description of the
resulting language. In this modification, the rightmost slot is
distinguished by not allowing either $\R$ or $\F$ code letters in
that slot. This ensures that there is always a slot present at the
right hand end -- an evolution may be complete when this is the only
remaining slot. With respect to this convention, the evolution of
$246153$ becomes:
\[
\slot  \to   \slot 1 \slot  \to  2 \slot 1 \slot  \to  2 \slot 1
\slot 3 \slot  \to  24 \slot 1 \slot 3 \slot  \to  24 \slot 153
\slot \to 246153 \slot.
\]
The corresponding encoding is $\M_1 \L_1 \M_2 \L_1 \F_2 \F_1$. For
the remainder of this paper, it is this variation of the insertion
encoding which we refer to as {\em the\/} insertion encoding.

We mention without proof the following result which will appear in
\cite{ALR:insertion}. It is not actually used in the next section,
but provides the motivation for it.

\begin{theorem}
Let $k$ be  a fixed positive integer. The collection of permutations
whose evolution requires at most $k$ slots at any point forms a
pattern class $\B_k$. The insertion encodings of $\B_k$ form a
regular language, as do the insertion encodings of any pattern class
$\B_k \intersect \avoid(X)$ where $X$ is a finite set of
permutations.
\end{theorem}

The theoretical methods of \cite{ALR:insertion} provide, in
principle, an effective method for determining the regular languages
representing the insertion encodings of $\B_k \intersect \avoid
(4231)$. In practice, these methods require various operations on
automata which are of exponential complexity and hence are
impractical for most values of $k$. 

Instead, in the next section, we describe a direct construction of
the automata which recognize words belonging to the insertion
encodings of elements of $\B_k \intersect \avoid (4231)$.

\section{The automata}

Consider a configuration of elements and slots which might arise in the evolution of a
$4231$ avoiding permutation. In this configuration there will be
some instances of patterns of the form $ \cdots \slot \cdots b
\cdots \slot  \cdots a \cdots$ where $b > a$. Wherever such an
instance occurs the first slot must be filled before the second slot
can be. Otherwise we would obtain four elements $\cdots d \cdots b
\cdots c \cdots a \cdots$ with $a < b < c < d$ in the resulting
permutation, that is, an instance of $4231$. Conversely, the only
way we could ever create a $4231$ pattern would be by insertion into
such a slot. Borrowing terminology from \cite{THEORY:queueJumping}
we say that in this configuration the second slot is {\em locked\/}
until such time as the first slot (and any other slots participating
in such patterns with it) are filled.

We now turn to the question of how locks are created, and how they
interact. Suppose that we have a configuration of $t$ slots:
\[
\alpha_1 \sslot \alpha_2 \sslot  \cdots \sslot \alpha_j \sslot
\alpha_{j+1} \sslot \cdots \alpha_t \sslot
\]
where $\alpha_1$ through $\alpha_t$ are certain sequences of
elements, $\alpha_1$ might be empty, but the remaining $\alpha$'s
are not. Suppose that the $j$th slot is not locked and we insert a
new maximum element $b$ into it, on the left for the sake of
argument. The new configuration is:
\[
\alpha_1 \sslot \alpha_2 \sslot  \cdots \sslot \alpha_j \, b \sslot
\alpha_{j+1} \sslot \cdots \alpha_t \sslot
\]
Taking  any slot from the first through the $(j-1)$st, $b$, any slot
from the $j$th through the $(t-1)$st and any element from $\alpha_t$
yields a $\slot b \slot a$ pattern. Thus all the slots from the
$j$th through the $(t-1)$st are now locked until all the slots from
the first through the $(j-1)$st have been filled.

We can record this in the new configuration by subscripting the
$j$th slot with the value $t-j$ -- which is to be read as ``the
$t-j$ consecutive slots beginning from this one are locked, until
the slots before it have been filled''. Alternatively, a more
attractive visual representation would be to place a bar over this
block of slots. Any slot under a bar cannot be filled, but bars are
removed when there are no slots to the left of them.

Other insertions into the $j$th slot create similar locks or bars.
If the intersection of two locks is non-empty then one must be
contained in the other, since a lock when created always begins at
the remaining slot just to the right of the current insertion and
ends at the penultimate slot.

It is possible for locks to be extended -- in the example above the
construction might proceed by adding a few more slots on the right
hand end (using middle insertions in the final slot), and then an
insertion on the right of the $(j-1)$st slot. Since this new lock
properly contains the old one, we can at this point discard the old
lock or simply extend its bar in the visual representation.

\begin{observation}
\label{OBS:locks}
If we know all the locking information about a configuration, then
we can determine which insertions are allowed. Furthermore, we can
determine the locking information of the configuration resulting
from any allowed insertion.
\end{observation}

The first part of the observation is trivial since, by definition,
insertions are allowed in the unlocked slots. The second follows
from the notes above, since the lock formed by any insertion does
not depend on the actual values present, only on the slots. Locks
are removed precisely when their left hand endpoint becomes the
leftmost slot.

By giving slots that are not at the left hand end of a lock a subscript of $0$ and then reading a
configuration only as a sequence of subscripts we see that the
configurations that can arise in the construction of a
$4231$-avoiding permutation are in one to one correspondence with
sequences $s_1 s_2 \cdots s_m$ (for $m \geq 1$) of non-negative
integers satisfying $s_1 = s_m = 0$, and if $s_k > 0$ then for all
$j < k + s_k$, $j + s_j \leq k + s_k$. The first condition expresses
the fact that the first and last slots are always unlocked, and the
second that if the $j$th slot lies within the lock on the $k$th
slot, then its lock cannot extend beyond the end of that one.
Sequences satisfying these conditions will be called {\em lock sequences\/}. It can
easily be established inductively (but is not actually required for
the following constructions) that every lock sequence can arise in
the evolution of some $4231$-avoiding permutation.

If we ignore the first and last slots (which can never be locked)
and think of the locks as subintervals of $\{ 1,2,3,\ldots,m\}$ we
see that they form a family of subintervals no two of which have the
same left endpoint, and with the property that if two intersect,
then one is a subinterval of the other. Of course this can be
thought of as a recursive description of how such arrangements of
locks can be created and it follows directly that the number of
configurations of locks on these $m$ elements is exactly the $m$th
large Schr\"{o}der number (sequence
\href{http://www.research.att.com/projects/OEIS?Anum=A006318}{A006318}
of \cite{EIS}). The large Schr\"{o}der numbers count paths in the
nonnegative half plane from $(0,0)$ to $(2n, 0)$ using steps
$\mathbf{u} = (1,1)$, $\mathbf{d} = (1,-1)$ and $\mathbf{h} =
(2,0)$. The correspondence is most easily seen from the set of such
paths to arrangements of locks. Associate the numbers $1$ through
$n$ with the $\mathbf{u}$'s and $\mathbf{h}$'s of such a sequence in
order. The locks are precisely the subintervals of numbers that
occur between some $\mathbf{u}$ and its matching $\mathbf{d}$. So,
for example the sequence $\mathbf{uhudhuhddh}$ corresponds to the
subintervals $[1,6]$, $[3,3]$ and $[5,6]$ of the interval $[1,7]$.

If we consider only locking sequences of length at most $k$ (for some fixed positive integer $k$) and the symbols of the insertion encoding which are allowed to operate on them, then Observation \ref{OBS:locks} and the discussion in the first paragraph of this section immediately imply the following result.

\begin{theorem}
Let $k$ be a fixed positive integer. There is a finite automation
$\Aut_k$ whose accepted language consists of the insertion encodings
of the permutations in  $\B_k \intersect \avoid (4231)$. The states
of $\Aut_k$ can be taken to be the lock sequences of length at most
$k$ and the transitions of $\Aut_k$ from a given sequence $s$ are
labelled by the codes of the allowed insertions in the slot
configuration corresponding to $s$, and are from $s$ to the lock
sequence labelling the result of the corresponding insertion.
\end{theorem}

The automata above are simply the restrictions of an automaton \Aut\
(with infinitely many states and an infinite language) that produces
the insertion encoding of all and only the elements of $\avoid
(4231)$. Its states are arbitrary lock sequences and its transitions
are precisely the allowed insertions within a lock sequence.

For illustrative purposes, consider $\Aut_4$. This  automaton has $10$ states represented by the lock sequences $0$, $00$, $000$, $010$, $0000$, $0010$, $0100$, $0110$, $0200$ and $0210$. A representative slot configuration for each of these states is: $\slot$, $\slot 1 \slot$, $\slot 1 \slot 2 \slot$, $\slot 2 \slot 1 \slot$, $\slot 1 \slot 2 \slot 3 \slot$, $\slot 1 \slot 3 \slot 2 \slot$, $\slot 2 \slot 1 \slot 3 \slot$, $\slot 2 \slot 1 4 \slot 3 \slot$, $\slot 3 \slot 1 \slot 2 \slot$ and $\slot 3 \slot 2 \slot 1 \slot$.  A complete transition table for this automaton is shown below. Each row illustrates the transitions available from the state specified at the left hand end of the row and double subscripts such as $\F_{12}$ indicate that both $\F_1$ and $\F_2$ induce the same transition.

{\small
\[
\begin{array}{r|cccccccccc}
 & 0 & 00 & 000 & 010 & 0000 & 0010 & 0100 & 0110 & 0200 & 0210 \\ \hline
0 & \L_1 & \M_1 \\
00 & \F_1 & \L_{12} \R_1 & \M_2 & \M_1 \\
000 & & \F_{12} & \L_{13} \R_2 & \R_1 \L_2  & \M_3 & \M_2 & & & \M_1 \\
010 & & \F_1 &  & \L_{13} \R_1 & & & \M_3 & & & \M_1 \\
0000 & & & \F_{13} & \F_2 & \L_{14} \R_3 & \R_2 \L_3 & & & \R_1 \L_2 \\
0010 & & & & \F_{12}  & & \L_{14} \R_2  & & & & \R_1 \L_2 \\
0100 & & & \F_1 & \F_3 & & & \L_{14} \R_3 & \L_3 & \R_1 \\
0110 & & & & \F_1 & & & & \L_{14} & & \R_1 \\
0200 & & & \F_1 & & & & & & \L_{14} \R_1 \\
0210 & & & & \F_1 & & & & & & \L_{14} \R_1
 \end{array}
\]
}

\section{Computational Methodology}

Let $\Lock$ be the set of all finite lock sequences. We order this set first by length, and then lexicographically within each length. This assigns an index (the position in this ordering) to each possible lock sequence. Armed with a table of Schr\"{o}der numbers, the recursive
description of $\Lock$ makes it relatively easy to compute these indices directly. Let 
\[
\ind : \Lock \to \mathbb{N}
\]
be the function which computes the index of a lock sequence.

Using $\ind$ and its inverse the states of \Aut\ can be indexed by the
natural numbers, and the transitions of \Aut\ can be determined. As
our goal is primarily to determine the growth rate of the language
accepted by $\Aut_k$ we can use these to construct the matrix $A_k$
whose entry in row $i$ and column $j$ is the number of transitions
between the state $\ind^{-1}(i)$ and $\ind^{-1}(j)$ (here $i$ and $j$ are
any pair of integers in the image of the lock sequences of length at
most $k$ under $\ind$).

The matrix $A_k$ is irreducible because the underlying directed multigraph is strongly connected.  Furthermore it is primitive as all the diagonal entries are non-zero (each state has a loop labelled  $\L_1$). Thus we can apply the Perron Frobenius theorem and conclude that $A_k$ has a unique dominant eigenvalue $\lambda_k$ which lies on the positive real axis and that the corresponding eigenvector is positive. Hence the limit 
\[
\lim_{n \to \infty} (e_1^T A_k^n e_1)^{1/n} = \lambda_k.
\]
Moreover, the generating
function for the language accepted by $\Aut_k$ is simply:
\[
\sum_{n = 0}^{\infty} e_1^T A_k^n e_1 t^n
\]
In other words, $\lambda_k$ is the growth rate of the language
accepted by $\Aut_k$ and hence the Wilf-Stanley limit of the class $\B_k \intersect \avoid
(4231)$.

The matrix $A_k$ is relatively sparse, so the eigenvalue $\lambda_k$
can be computed without great difficulty even for moderately large
values of $k$. For instance, if $k = 13$, $A_k$ is a square matrix
with $6589728$ rows. There are at most $46$ transitions from any
state in the automaton (this is achieved in the state with 12 slots
having no locks -- many states have significantly fewer
transitions). However, no row has quite this many non zero entries
as there are always several transitions to the same state. 

\begin{theorem}
The Wilf-Stanley limit, $\ws{4231}$ is at least $9.35$.
\end{theorem}

\begin{pf}
Let $A = A_{13}$. Because $A$ is irreducible, primitive, and non-negative an iterative scheme to compute its dominant eigenvalue is guaranteed to converge. That is, we may define a sequence of vectors $\vec{v}_k$ where $\vec{v}_1 = e_1$ and $\vec{v}_{k+1}$ is a scalar multiple of $A \vec{v}_{k}$  having some fixed norm. This method, implemented in Java, produced a dominant eigenvalue of $9.3508$ for the matrix $A$ together with an approximate eigenvector $\vec{v}$. Direct computation then showed that:
$
A \vec{v} \geq (9.35) \vec{v}.
$
Since the entries of $A$ are all non-negative and the diagonal entries are all positive it follows that 
$
A^n \vec{v} \geq (9.35)^n \vec{v}
$
for all positive integers $n$. Since the first coordinate of $\vec{v}$ is non-zero, it also follows that:
\[
\lim_{n \to \infty} (e_1^T A_k^n e_1)^{1/n}  \geq 9.35
\]
which, as noted above, establishes the claim of the theorem. \qed
\end{pf}

The values of $s_n(4231)$ are reported for $n \leq 20$ as sequence \href{http://www.research.att.com/projects/OEIS?Anum=A061552}{A061552} in \cite{EIS}. The recursive method used to compute these numbers is described in \cite{Marinov:EJC2002} and its exact complexity has not been analysed. A permutation requiring more than $k$ slots to produce in the insertion encoding must have length at least $2k$ so $s_n(1324)$ is the $(1,1)$ entry of $A_k^n$ for any $k > n/2$. Choosing $k = 13$ allows us to report that the values of the sequence $s_n(4231)$ for $n$ between $21$ and $25$ are:
$1535346218316422$, 
$12015325816028313$, 
$94944352095728825$,
$757046484552152932$ and 
$6087537591051072864$.

As $A_k$ has asymptotically $O(k(1+\sqrt{2})^{2k})$ non-zero entries, the complexity of the computation of $s_n(4231)$ by this method is not more than $O(n^2 (1 + \sqrt{2})^n)$ (and the constants are not large).

\section{Conclusions}

The lower bounds presented here leave the question of the true
growth rate of $\avoid (4231)$ intriguingly open. Since the sequence $\lambda_k$ is monotone increasing and bounded above by $\ws{4231}$ it has a limit $\lambda_{\infty} \leq \ws{4231}$. Although the generating functions for the language accepted
by $\Aut_k$ and $\avoid (4231)$ agree through at the first
$2k$ terms, this does not necessarily guarantee that $\lambda_{\infty} = \ws{4231}$. So this raises:

\begin{question}
\label{Q:Limit}
Is $\lim_{k \to \infty} \lambda_k = \ws{4231}$?
\end{question}

The growth rates
of the automata languages for different values of $k$ are presented
below.
\[
\begin{array}{cc}
k & \lambda_k \\ \hline 
1 & 1.0000 \\ 
2 & 3.4142 \\
3 & 5.1120 \\
4 & 6.2262 \\
5 & 7.0014 \\
6 & 7.5693 \\
7 & 8.0029 \\
8 & 8.3450 \\
9 & 8.6220 \\
10 & 8.8511 \\
11 & 9.0439 \\
12 & 9.2085 \\
13 & 9.3508
\end{array}
\]
We leave it to the reader to decide how to extrapolate this sequence.
However, the value obtained will depend on how one models the behaviour of
the difference $\lambda_{\infty} - \lambda_k$ as a function of $k$. Our best guess, based on an empirical observation that the plot of $1/\sqrt{k}$ against $\lambda_k$ is roughly linear is that the limiting value lies between $11$ and $12$. Computing $\lambda_k$
for larger values of $k$ is possible -- though we note that there
are  over $6.5 \times 10^6$ states in $\Aut_{13}$ and the number of states goes up by a factor of roughly $5.8$ for each additional slot so significant
further progress in this direction is limited by the obvious
combinatorial explosion. However, the natural structure of the
states of \Aut\ leaves open the possibility of a closed form or
asymptotic analysis of the limiting case. 

We suspect that the answer to Question \ref{Q:Limit} is yes, but the evidence is not entirely convincing. It consists of the observation that for  $\avoid(312)$ we can carry out a similar analysis (and of course we know that $\ws{312} = 4$) and because of the simple form of the corresponding automaton which only has one state for each number of slots, we can prove that the maximal eigenvalues do converge to $4$. On the other hand for the class $\avoid(4321)$ with $\ws{4321} = 9$ it is also the case that the underlying automata are relatively simple, the one for $k$ slots having only $O(k^2)$ states. The corresponding dominant eigenvalues do appear to converge to $9$ but the rate of convergence is quite slow.

The results above show that the class of $4231$ avoiders has
strictly larger growth rate than any other class avoiding a single
permutation of length $4$. This throws open once again the question of what makes one pattern harder to avoid than another. That is:

\begin{question}
Among the classes $\avoid(\pi)$ where $\pi$ is a single permutation
of length $k$, which have the largest growth rates? What is this largest growth rate? More generally, given two permutations $\pi$ and $\tau$ are there general methods for deciding whether or not $\ws{\pi} \geq \ws{\tau}$?
\end{question}

\bibliographystyle{elsart-num}
\bibliography{av4231}

\end{document}